\documentclass[12pt,twoside]{amsart}    

\usepackage{amssymb}

\def\qed{$\rlap{$\sqcap$}\sqcup$}
\usepackage{latexsym}

\begin{document}           

\begin{center}
{\ }\\
{\Large {\bf On the degree two entry of a Gorenstein $h$-vector and
a  conjecture of Stanley}} \\ [.250in]
{JUAN MIGLIORE\\
Department of Mathematics, University of Notre Dame, Notre Dame, IN
46556, USA
\\E-mail: Juan.C.Migliore.1@nd.edu}
{{\ }
\\UWE NAGEL\footnote[1]{The second author gratefully acknowledges partial
support from
and the hospitality of the Institute for Mathematics \& its
Applications
at the University of Minnesota.}\\
Department of Mathematics, University of Kentucky, 715 Patterson
Office Tower,\\Lexington, KY 40506-0027, USA
\\E-mail: uwenagel@ms.uky.edu}
{{\ }
\\FABRIZIO ZANELLO\\
Department of Mathematics, University of Notre Dame, Notre Dame,  IN
46556, USA
\\E-mail: zanello@math.kth.se\\
Current address: Department of Mathematical Sciences,\\Michigan
Technological University, Houghton, MI 49931-1295, USA}
\end{center}
{\ }\\
\\{\small
ABSTRACT. In this short note we establish a (non-trivial) lower
bound on the degree two entry $h_2$ of a Gorenstein $h$-vector of
any given socle degree $e$ and any codimension $r$.\\
In particular, when $e=4$, that is for Gorenstein $h$-vectors of the
form $h=(1,r,h_2,r,1)$, our lower bound allows us to prove a conjecture
of Stanley on the order of magnitude of the minimum value, say $f(r)$,
that $h_2$ may assume. In fact, we show that
$$
\lim_{r\rightarrow \infty} \frac{f(r)}{ r^{2/3}}= 6^{2/3}.$$
In general, we wonder whether our lower bound is sharp for all integers
$e\geq 4$ and $r\geq 2$.
} \\


{\ }\\

\section{Introduction}

Gorenstein rings arise in many areas of mathematics - including algebraic geometry, combinatorics, and complexity theory (see, e.g, $[Hu]$, $[Pa]$, $[St4]$, $[At]$, $[St5]$, $[KS]$). Often these rings are standard graded algebras over a field $k$. Then it is a basic problem to understand the vector space dimensions of the graded components of such  algebras, i.e.\ their Hilbert functions. While there is a complete classification of the Hilbert functions of Cohen-Macaulay algebras due to Macaulay, there is not even a conjecture about the possible Hilbert functions of Gorenstein algebras. Hence, it is a natural first step to ask for bounds on the Hilbert function in terms of some information on the given algebra.
This note addresses some instances of this problem. Our results are strong enough to settle a conjecture of Stanley that will be discussed below.


The power series generated by the Hilbert function of a standard graded Gorenstein $k$-algebra $A$ provides the so-called {\em $h$-vector} $h=(h_0,\ldots,h_e)$ of $A$. This finite tuple of positive integers, together with the (Krull) dimension of $A$, completely describes  the Hilbert function of $A$. Equivalently, we may assume without loss of generality that $A$ is artinian. Then the $h$-vector of $A$ consists of the positive integers $h_i := \dim_k A_i$. By duality, it is symmetric about the middle, and in particular $h_e = h_0 = 1$. We call $h_1$ the {\em codimension} of $A$, and $e$ the {\em socle degree} of $A$.

In codimension $h_1 \leq 3$, all Gorenstein $h$-vectors have been characterized (see $[St1]$ and also $[Za1]$), and it turned out that they are exactly the so-called {\em SI-sequences} (named after Stanley and Iarrobino). In particular, they are {\em unimodal}, i.e. they never strictly increase after a strict decrease. While it is open whether non-unimodal Gorenstein $h$-vectors of codimension 4 exist (though see $[IS]$ for a recent unimodality result), it is known that non-unimodal Gorenstein $h$-vectors exist in every codimension $h_1 \geq 5$ (see $[BI]$, $[BL]$, $[Bo]$). The first such example was given by Stanley (see $[St1]$, Example 4.3). He showed that $h=(1,13,12,13,1)$ is indeed a Gorenstein $h$-vector. Non-unimodality occurs here in degree two. This note focuses on the degree two entry of Gorenstein $h$-vectors.

It is well-known that  the degree 2 entry satisfies $h_2 \leq \binom{h_1 + 1}{2}$. However,  if  the socle degree $e$ is at least 4, very little is known on the least value that $h_2$ may assume. In particular, there is not even a precise conjecture about it.
The first main result of this note (Theorem 4) gives a non-trivial lower bound on $h_2$ in terms of the codimension $h_1$ and the socle degree $e$. Its proof relies on  results of Macaulay and Green as well as
on an observation of Stanley.
Considerably extending results of Stanley in $[St3]$, our lower bound immediately implies that certain Gorenstein algebras of small codimension and small socle degree always have  unimodal $h$-vectors (Corollary 6).

In Section 3 we apply our bound to investigate Gorenstein algebras of
 socle degree 4; that is, we consider Gorenstein $h$-vectors of the form $h=(1,r,h_2,r,1)$.    Fixing the codimension $r$, denote the least possible value that $h_2$ may assume by $f(r)$. Stanley ($[St2]$) conjectured that the limit ${\displaystyle \lim_{r\rightarrow \infty} \frac{f(r)}{r^{2/3}}}$ exists and he also conjectured its precise value. However, the existence of the limit was open and, assuming its existence, bounds for it were shown by Stanley ($[St3]$) and Kleinschmidt ($[Kl]$). In Theorem 9 we  prove Stanley's conjecture by showing that
$$
\lim_{r\rightarrow \infty} \frac{f(r)}{ r^{2/3}}= 6^{2/3}.
$$
It follows that our lower bound is asymptotically optimal in this case.

\section{The lower bound}

We begin by recalling  results of Macaulay, Green, and Stanley that we will need in this paper.\\
\\\indent
{\bf Definition-Remark 1.} Let $n$ and $i$ be positive integers. The {\em i-binomial expansion of n} is $$n_{(i)}=\binom{n_i}{ i}+\binom{n_{i-1}}{ i-1}+...+\binom{n_j}{ j},$$ where $n_i>n_{i-1}>...>n_j\geq j\geq 1$. Such an expansion always exists and it is unique (see, e.g.,  $[BH]$, Lemma 4.2.6).\\\indent
Following $[BG]$, we define, for any integers $a$ and $b$,
$$
(n_{(i)})_{a}^{b}=\binom{n_i+b}{ i+a}+\binom{n_{i-1}+b}{ i-1+a}+...+\binom{n_j+b}{ j+a},
$$
where we set $\binom{m}{ q}=0$ whenever $m<q$ or $q<0$.\\

Each standard graded $k$-algebra $A$ can be written as $A = R/I$ where $R = k[x_1,\ldots,x_r]$ is some polynomial ring and $I \subset R$ is  a homogeneous ideal. Since we are interested in the Hilbert function of $A$, we may and will assume that the field $k$ is infinite.
\\

\indent
{\bf Theorem 2} {\em Let $L \in A$ be a general linear form. Denote by $h_d$  the degree $d$ entry of the Hilbert function of $A$ and by $h_d^{'}$ the degree $d$ entry of the Hilbert function of $A/L A$. Then:\\\indent
i) {\em (Macaulay)} $$h_{d+1}\leq ((h_d)_{(d)})_{1}^{1}.$$\indent
ii) {\em (Green)} $$h_d^{'}\leq ((h_d)_{(d)})_{0}^{-1}.$$}
\\\indent
{\bf Proof.} i) See $[BH]$, Theorem 4.2.10.\\\indent
ii) See $[Gr]$, Theorem 1.{\ }{\ }\qed \\
\\\indent
The following  simple observation is not new (see for instance $[St3]$, bottom of p. 67):\\
\\\indent
{\bf Lemma 3} (Stanley). {\em Let $A$ be an artinian Gorenstein algebra, and let $L\in A$ be any linear form. Then the  $h$-vector of $A$ can be written as
$$h := (h_0,\ldots,h_e) = (1, b_1+c_1,..., b_e+c_e = 1),$$ where
$$b=(b_1=1,b_2,...,b_{e-1},b_e=1)$$
is the $h$-vector of $A/(0:L)$ (with the indices shifted by 1 to the left), which is a Gorenstein algebra, and
$$c=(c_0=1,c_1,...,c_{e} = 0)$$
is the $h$-vector of $A/L A$.}\\
\\\indent
{\bf Proof.} The fact that $A/(0:L)$ is Gorenstein of socle degree $e-1$ is well-known and follows, for example, by computing its socle. The decomposition of $h$ is an immediate consequence of the exact sequence
$$0\longrightarrow A/(0:L) (-1)\stackrel{\cdot L}{\longrightarrow} A \stackrel{}{\longrightarrow} A/ L A \longrightarrow 0$$
that is induced by the multiplication by $L$.{\ }{\ }\qed \\
\\\indent
For all integers $r\geq 2$ and $e\geq 4$, we define $f(r,e)$ as the least possible value of the degree two entry $h_2$ among the Gorenstein $h$-vectors $h = (1, r, h_2,\ldots,h_e = 1)$ of codimension $r$ and socle degree $e$. The first main result of this note is the following estimate:\\
\\\indent
{\bf Theorem 4.} {\em $$f(r,e)\geq (r_{(e-1)})_{-1}^{-1}+(r_{(e-1)})_{-(e-3)}^{-(e-2)}.$$}
\\\indent
{\bf Proof.}
Let $h = (1,r,h_2,\ldots,h_e = 1)$ be a  Gorenstein vector such that $h_2 = f(r,e)$. Let $L \in R$ be a general linear form.
By Lemma 3, we can decompose $h$ as $$b+c=(0,b_1=1,b_2,...,b_{e-1},b_e)+(1,c_1,c_2,...,c_{e-1},c_e),$$ where $b$ is a Gorenstein $h$-vector of socle degree $e-1$ (shifted by 1 to the right), and $c$ is the $h$-vector of some artinian algebra. In particular, $b_2=b_{e-1}$. By Green's Theorem 2, ii), we have $c_{e-1}\leq (r_{(e-1)})_{0}^{-1}$. It follows that
$$
b_2=b_{e-1}\geq r-(r_{(e-1)})_{0}^{-1}=(r_{(e-1)})_{-1}^{-1}.
$$
Hence we get $c_2\leq f(r,e)-(r_{(e-1)})_{-1}^{-1}$.  Using repeatedly  Macaulay's Theorem 2, i) and the fact that $(m_{(d)})_1^1$ is an increasing function in $m$, we can compare $c_{e-1}$ and $c_2$, and we obtain
\begin{equation}\label{mm}c_{e-1}\leq ((f(r,e)-(r_{(e-1)})_{-1}^{-1})_{(2)})^{e-3}_{e-3}.\end{equation}\indent
We now claim that $g(r,e):=\left((f(r,e)-(r_{(e-1)})_{-1}^{-1})_{(2)}\right)^{e-3}_{e-3} \geq (r_{(e-1)})_{0}^{-1}$.

Assume the contrary, i.e. that $g(r,e)< (r_{(e-1)})_{0}^{-1}$. Then, by (\ref{mm}), $c_{e-1}$ is less than $(r_{(e-1)})_{0}^{-1}$, say $c_{e-1}=(r_{(e-1)})_{0}^{-1}-a$, for some integer $a>0$. Hence,
$$b_2=b_{e-1}=r-c_{e-1}=r-(r_{(e-1)})_{0}^{-1}+a=(r_{(e-1)})_{-1}^{-1}+a,$$
and therefore $c_2=f(r,e)-(r_{(e-1)})_{-1}^{-1}-a$. Repeated use of Macaulay's Theorem 2, i) provides
$$((f(r,e)-(r_{(e-1)})_{-1}^{-1}-a)_{(2)})^{e-3}_{e-3}=((c_2)_{(2)})^{e-3}_{e-3}\geq c_{e-1}=(r_{({e-1})})_{0}^{-1}-a.
$$
Since, for any given $d$, $(m_{(d)})^1_1$ is a actually a strictly increasing function in $m$, it follows that $$g(r,e)=((f(r,e)-(r_{(e-1)})_{-1}^{-1})_{(2)})^{e-3}_{e-3} \geq ((f(r,e)-(r_{(e-1)})_{-1}^{-1}-a)_{(2)})^{e-3}_{e-3} + a
\geq (r_{(e-1)})_{0}^{-1},
$$
a contradiction. This proves the claim.\\

Using again monotonicity, our claim $((f(r, e)-(r_{(e-1)})_{-1}^{-1})_{(2)})^{e-3}_{e-3}\geq (r_{(e-1)})_{0}^{-1}$ implies
$$f(r,e)-(r_{(e-1)})_{-1}^{-1}\geq (r_{(e-1)})^{-(e-2)}_{-(e-3)},$$
and the theorem follows.{\ }{\ }\qed \\
\\\indent
We now record some interesting consequences of Theorem 4.\\
\\\indent
{\bf Example 5.} Consider the case $r=13$ and $e=4$. Our Theorem 4 yields $f(13,4)\geq 12$. Hence, Stanley's Example 4.3 in $[St1]$, $(1,13,12,13,1)$, shows that $f(13, 4) = 12$. In particular, our bound is sharp in this case. \\

Our lower bound allows us to show that certain Gorenstein $h$-vectors must be unimodal. Part i) below is a considerable extension of Fact (c) on page 67 in $[St3]$.
\\

{\bf Corollary 6.} {\em i).  If $r\leq 9$,  then, for all $e\geq 4$, $f(r,e)=r$. In particular, every socle degree 4 Gorenstein $h$-vector of codimension $r \leq 9$ is unimodal.\\\indent
ii). For all $r\leq 13$, $f(r,5)=r$. In particular, every socle degree 5 Gorenstein $h$-vector of codimension $r \leq 13$ is unimodal.\\}
\\\indent
{\bf Proof.} In both cases a
 standard computation shows that Theorem 4 provides $f(r,e)\geq r$. Since $(1,r,r,\ldots,r,1)$ is a Gorenstein $h$-vector whenever $r \geq 1$, we conclude $f(r,e)=r$, as desired.{\ }{\ }\qed \\
\\\indent
{\bf Remark 7.}  One might think that, fixing the codimension $r$, it becomes easier to find non-unimodal Gorenstein $h$-vectors, the larger the socle degree is. However, this is not always true (at least for short $h$-vectors). Indeed, Example 5 and Corollary 6, ii) show that in case $r = 13$, the existence of a non-unimodal Gorenstein $h$-vector of socle degree $4$ does not imply the existence of non-unimodal Gorenstein $h$-vectors of socle degree 5.
\\

However, fixing the socle degree, but allowing the codimension to increase, non-unimo\-dality can be carried over. More precisely, we have: \\

{\bf Proposition 8.} {\em Given any integer $e\geq 4$, the function $r-f(r,e)$ is non-decreasing in~ $r$.\\\indent
In particular, since $(1,13,12,13,1)$ is a non-unimodal Gorenstein $h$-vector, there exist non-unimodal Gorenstein $h$-vectors of codimension $r$ and socle degree 4 for all $r\geq 13$.}\\
\\\indent
{\bf Proof.}
The result is a consequence of the following claim:\\
{\em Claim.} If  $(1,r,h_2,...,r,1)$ is a Gorenstein $h$-vector, then $(1,r+1,h_2+1,...,r+1,1)$ is also a Gorenstein $h$-vector.

Indeed, assuming the claim, we immediately have that $f(r+1,e)\leq f(r,e)+1$, and therefore  $$r+1-f(r+1,e)\geq r+1-[f(r,e)+1]=r-f(r,e),$$ as desired.

{\em Proof of Claim.} Let $(1,r,h_2,...,r,1)$ be the $h$-vector of the Gorenstein algebra  $A = R/I$. Using Macaulay's inverse systems if the base field $k$ has characteristic 0, and divided powers if $k$ has positive characteristic (see, e.g., $[Ge]$ or $[IK]$), it is well-known that $I$ is the annihilator of a homogeneous  polynomial $F \in k[y_1,y_2,...,y_r]$ of degree $e$.

Now consider the form $\bar{F} := F+y_{r+1}^e$ in $r+1$ variables, and let $I' \subset R' :=k[x_1,\ldots,x_{r+1}]$ be the annihilator of $\bar{F}$. Note that a form $g \in R'$ annihilates $F$ if and only if $g,x_{r+1} g\in I'$, unless $g = x_{r+1}^j$ for some $j \leq e$.  It easily follows that the $h$-vector of $R'/I'$ is
 $(1,r+1,h_2+1,...,r+1,1)$, as claimed.{\ }{\ }\qed

\section{Algebras of socle degree 4}

We now turn to an application of Theorem 4 to Gorenstein algebras of socle degree 4, i.e.\ to Gorenstein $h$-vectors of the form $(1,r,h_2,r,1)$. In this case, the lower bound of Theorem 4 becomes
$$
h_2 \geq f(r):=f(r,4)\geq (r_{(3)})_{-1}^{-1}+(r_{(3)})_{-1}^{-2}.
$$
On page 79 in  $[St2]$, Stanley conjectured that for $r \to \infty$ the function $\frac{f(r)}{r}$ converges and that the limit $c$ is $6^{2/3}$. Assuming that the limit $c$ exists, Stanley himself showed that $\frac{1}{ 2}6^{2/3}\leq c\leq \frac{3}{2} 6^{2/3}$ (see $[St3]$, p.\ 68).  The upper bound has been improved by Kleinschmidt who showed $c \leq 6^{2/3}$ (see $[Kl]$), provided the limit $c$ exists.

The goal of this section is to show Stanley's conjecture by
using Theorem 4 for $e=4$, along with a suitable refinement of Stanley's use of trivial extensions:\\
\\\indent
{\bf Theorem 9.} {\em $$\lim_{r\rightarrow \infty} \frac{f(r)}{ r^{2/3}}=6^{2/3}.$$}
\\\indent
First we need two preparations. The first, once again due to Stanley, determines the $h$-vector of the trivial extension of a level algebra. Recall that an artinian $k$-algebra is called {\em level} if its socle is concentrated in one degree. \\
\\\indent
{\bf Lemma 10.} (Stanley). {\em Given a level algebra with $h$-vector $h=(1,h_1,...,h_j)$, there exists a Gorenstein algebra (called its {\em trivial extension}) having the $h$-vector\\ $H=(1,H_1,...,H_j,H_{j+1}=1)$, where, for $i=1,...,j$, $$H_i=h_i+h_{j+1-i}.$$}
\\\indent
{\bf Proof.} See $[St1]$, Example 4.3, and also $[Re]$. For a proof using inverse systems, see $[BI]$, Lemma 1.{\ }{\ }\qed \\
\\\indent
{\bf Lemma 11.} {\em Every integer $r \geq 4$ can be written as
$$
r = m+\binom{m+1}{ 3}+\binom{a+1}{ 2}+b,
$$
where $m, a, b$ are integers satisfying $1\leq a\leq  m-1$, and $0\leq b\leq a+2$.}\\
\\\indent
{\bf Proof.} Let $m$ be the largest integer such that $m + \binom{m+1}{3} \leq r$. If $m + \binom{m+1}{3} < r$, then we read off $a$ and $b$ from the 2-binomial expansion of $r - m - \binom{m+1}{3}$. Otherwise, we get the claimed presentation from the identity
$$
r = m + \binom{m+1}{3} = m-1 +  \binom{m}{3} + \binom{m-1}{2} + m. \quad \mbox{\qed}
$$
 \\

We are ready to establish our second main result.
\\

{\bf Proof of Theorem 9.} We want to show that $$F(r):=\frac{f(r)}{ (6r)^{2/3}}$$ converges to 1 as $r$ goes to infinity. To this end we will  exhibit functions $G$ and $H$ such that, for all $r$, $G(r)\leq F(r)\leq H(r)$ and both $G$ and $H$ converge to 1.

First, we use Theorem 4 to obtain the lower bound.
The 3-binomial expansion of $r$ can be rewritten as $r =\binom{k}{ 3}+\binom{p}{ 2}+\binom{q}{ 1}$, where $k > p, q \geq 0$. In particular, we have $r \leq \binom{k}{3} + \binom{k-1}{2} + \binom{k-1}{1}$. Furthermore, Theorem 4 with $e = 4$ provides
$$
f(r) \geq (r_{(3)})_{-1}^{-1}+(r_{(3)})_{-1}^{-2} \geq \binom{k-1}{2} + \binom{k - 2}{2} = k^2 - 4k + 4.
$$
It follows that
$$
F(r) \geq G(r) := \frac{k^2 - 4k + 4}{6^{2/3} \cdot \left [\binom{k}{3} + \binom{k-1}{2} + \binom{k-1}{1} \right ]^{2/3}}.
$$
Since $k \to \infty$ as $r \to \infty$, we conclude that ${\displaystyle \lim_{r\rightarrow \infty} G(r) = 1}$.

In order to establish the desired upper bound, we will construct suitable Gorenstein algebras.
Let $r \geq 4$. By Lemma 11, there are integers $m, a, b$ such that $1\leq a\leq  m-1$,  $0\leq b\leq a+2$, and $r = m+\binom{m+1}{ 3}+\binom{a+1}{ 2}+b$. By Macaulay's Theorem 1, i), there is a $k$-algebra $A$ with $h$-vector
$$
h=(1, m, \binom{m}{ 2}+a, \binom{m+1}{ 3}+\binom{a+1}{ 2}).
$$
Note that every such algebra $A$ must be a level algebra (see $[Za2]$, Theorem 3.5, or $[Na]$, Corollary 3.9).
Hence, by Lemma 10, there exists a Gorenstein $h$-vector of the form:
$$
(1, m+\binom{m+1}{ 3}+\binom{a+1}{ 2},2\binom{m}{ 2}+2a,m+\binom{m+1}{ 3}+\binom{a+1}{ 2},1).
$$
Using the construction in the proof of Proposition 8,  we see that there is a Gorenstein $h$-vector
$$
(1, r = m+\binom{m+1}{ 3}+\binom{a+1}{ 2}+b, 2\binom{m}{ 2}+2a+b, r, 1).
$$
It follows that
$$
f(r) \leq 2\binom{m}{ 2}+2a+b  \leq 2 \binom{m}{2} + 2 (m-1) + m + 1 = m^2 + 2 m -1. $$
Using $r \geq m + \binom{m}{3}$, we obtain
$$ F(r) \leq H(r) := \frac{m^2 + 2 m -1}{6^{2/3} \left [ m + \binom{m}{3} \right ]^{2/3}}.
$$
Since $m \to \infty$ as $r \to \infty$, we get ${\displaystyle \lim_{r\rightarrow \infty} H(r) = 1}$.
This  completes the proof of the theorem.{\ }{\ }\qed \\
\\\indent
The proof of the above theorem shows that the lower bound in Theorem 4 in case  $e = 4$ and $f(r) = f(r, 4)$ agree for $r \to \infty$ up to terms of lower order. Thus, we conclude this note by asking:\\
\\\indent
{\bf Question 12.} For which values of $r$ and $e$ is the lower bound on $f(r,e)$ supplied by Theorem 4 sharp?\\

In particular, the first interesting case would be to determine if the sequence $h =
(1,10,9,10,1)$ is a Gorenstein $h$-vector. Note that since 9 is odd, Stanley's
trivial extension method (cf.\ Lemma 10) does not work.\\
\\
{\bf Remark added after acceptance of paper.} In a work currently in progress, M. Boij and the third author have recently proved that several new $h$-vectors, including $(1,10,9,10,1)$, are {\em not} Gorenstein, and hence that the lower bound of Theorem 4 is in general not sharp.\\
\\
\\
\begin{center}
{REFERENCES}
\end{center}
\bigskip

{\small
\noindent
$[At]$ {\ } C.A. Athanasiadis: {\it Ehrhart polynomials, simplicial polytopes, magic squares and a conjecture of Stanley},
J. Reine Angew. Math. {\bf  583} (2005), 163-174.\\
$[BI]$ {\ } D. Bernstein and A. Iarrobino: {\it A nonunimodal graded Gorenstein Artin algebra in codimension five}, Comm.  Algebra {\bf 20} (1992), No. 8, 2323-2336.\\
$[BG]$ {\ } A.M. Bigatti and A.V. Geramita: {\it Level Algebras, Lex Segments and Minimal Hilbert Functions}, Comm. Algebra {\bf 31} (2003), 1427-1451.\\
$[Bo]$ {\ } M. Boij: {\it Graded Gorenstein Artin algebras whose Hilbert functions have a large number of valleys}, Comm. Algebra {\bf 23} (1995), No. 1, 97-103.\\
$[BL]$ {\ } M. Boij and D. Laksov: {\it Nonunimodality of graded Gorenstein Artin algebras}, Proc. Amer. Math. Soc. {\bf 120} (1994), 1083-1092.\\
$[BH]$ {\ } W. Bruns and J. Herzog: {\it Cohen-Macaulay rings}, Cambridge studies in advanced mathematics {\bf 39}, Revised edition (1998), Cambridge, U.K..\\
$[Ge]$ {\ } A.V. Geramita: {\it Inverse Systems of Fat Points: Waring's Problem, Secant Varieties and Veronese Varieties and Parametric Spaces of Gorenstein Ideals}, Queen's Papers in Pure and Applied Mathematics {\bf 102}, The Curves Seminar at Queen's (1996), Vol. X, 3-114.\\
$[Gr]$ {\ } M. Green: {\it Restrictions of linear series to hyperplanes, and some results of Macaulay and Gotzmann}, Algebraic curves and projective geometry (1988), 76-86, Trento; Lecture Notes in Math. {\bf 1389} (1989), Springer, Berlin.\\
$[Hu]$ {\ } C. Huneke: {\it Hyman Bass and Ubiquity: Gorenstein Rings}, Contemp. Math. {\bf 243} (1999), 55-78.\\
$[IK]$ {\ } A. Iarrobino and V. Kanev: {\it Power Sums, Gorenstein Algebras, and Determinantal Loci}, Springer Lecture Notes in Mathematics {\bf 1721} (1999), Springer, Heidelberg.\\
$[IS]$ {\ } A. Iarrobino and H. Srinivasan: {\it Some Gorenstein Artin algebras of embedding dimension four, I: components of $PGOR(H)$ for $H=(1,4,7,...,1)$}, J. Pure Appl. Algebra {\bf 201} (2005), 62-96.\\
$[Kl]$ {\ } P.\ Kleinschmidt: {\it \"Uber Hilbert-Funktionen graduierter Gorenstein-Algebren}, Arch.\ Math. {\bf 43} (1984), 501-506. \\
$[KS]$ {\ } A.R. Klivans and A. Shpilka: {\it Learning arithmetic circuits via partial derivatives}, in: Proc. 16th Annual Conference on Computational Learning Theory, Morgan Kaufmann Publishers (2003), 463-476.\\
$[Na]$ {\ } U.\ Nagel: {\em Empty simplices of polytopes and graded
  Betti numbers}, Discrete Comput.\ Geom.\ (to appear). \\
$[Pa]$ {\ } R. Pandharipande: {\it Three questions in Gromov-Witten theory}, in: Proceedings of the International Congress of Mathematics, Vol. II, Higher Ed. Press, Beijing (2002), 503-512.\\
$[Re]$ {\ } I. Reiten: {\it The converse to a theorem of Sharp on Gorenstein modules}, Proc. Amer. Math. Soc. {\bf 32} (1972), 417-420.\\
$[St1]$ {\ } R. Stanley: {\it Hilbert functions of graded algebras}, Adv. Math. {\bf 28} (1978), 57-83.\\
$[St2]$ {\ } R. Stanley: {\it Combinatorics and Commutative Algebra}, First Ed., Progress in Mathematics {\bf 41} (1983), Birkh\"auser, Boston, MA.\\
$[St3]$ {\ } R. Stanley: {\it Combinatorics and Commutative Algebra}, Second Ed., Progress in Mathematics {\bf 41} (1996), Birkh\"auser, Boston, MA.\\
$[St4]$ {\ } R. Stanley: {\it A monotonicity property of $h$-vectors and $h^{*}$-vectors}, European J. Combin.  {\bf 14} (1993), 251-258.\\
$[St5]$ {\ } R. Stanley: {\it The number of faces of a simplicial convex polytope}, Adv. Math. {\bf 35} (1980), 236-238.\\
$[Za1]$ {\ } F. Zanello: {\it Stanley's theorem on codimension 3 Gorenstein $h$-vectors}, Proc. Amer. Math. Soc. {\bf 134}
(2006), No. 1, 5-8.\\
$[Za2]$ {\ } F. Zanello: {\it When is there a unique socle-vector associated to a given $h$-vector?}, Comm. in Algebra {\bf 34} (2006), No. 5, 1847-1860.

}

\end{document}